\documentclass[fleqn,10pt]{wlscirep}
\usepackage[utf8]{inputenc}
\usepackage[T1]{fontenc}

\title{Homological Scaffold via Minimal Homology Bases}

\usepackage{amsthm,amsmath, amssymb}
\usepackage{mathtools}
\usepackage{placeins}
\usepackage{verbatim}
\usepackage{amsthm}
\usepackage{bbm}
\usepackage{newtxmath}
\usepackage{todonotes}
 
\theoremstyle{definition}
\newtheorem{definition}{Definition}[section]
\DeclareMathOperator{\im}{Im}

\author[1,*]{Marco Guerra}
\author[1]{Alessandro De Gregorio}
\author[3]{Ulderico Fugacci}
\author[2]{Giovanni Petri}
\author[1]{Francesco Vaccarino}
\affil[1]{Politecnico di Torino, Department of Mathematical Sciences, Torino, Italy}
\affil[2]{ISI Foundation, Torino, Italy}
\affil[3]{CNR-IMATI, Genova, Italy}

\affil[*]{marco.guerra@polito.it}


\begin{abstract}
The homological scaffold leverages persistent homology to construct a topologically sound summary of a weighted network. However, its crucial dependency on the choice of representative cycles hinders the ability to trace back global features onto individual network components, unless one provides a principled way to make such a choice. In this paper, we apply recent advances in the computation of minimal homology bases to introduce a quasi-canonical version of the scaffold, called minimal, and employ it to analyze data both real and in silico. At the same time, we verify that, statistically, the standard scaffold is a good proxy of the minimal one for sufficiently complex networks. 
\end{abstract}
\begin{document}

\flushbottom
\maketitle
\thispagestyle{empty}

%\noindent Please note: Abbreviations should be introduced at the first mention in the main text – no abbreviations lists. Suggested structure of main text (not enforced) is provided below.

\section{Introduction}

\noindent
Network science has long represented the cornerstone theory in dealing with complex, heterogeneous multi-agent systems.
Network descriptions have found wide applications and had a significant impact on a wide range of fields (\cite{newman2003structure,barrat2004architecture}), including social networks (\cite{granovetter1977strength,vega2007complex}), epidemiology (\cite{pastor2015epidemic,colizza2006role}), biology (\cite{girvan2002community,alon2003biological}), and neuroscience (\cite{bassett2017network,bullmore2009complex,bassett2006small}). 

In recent years, new approaches to the analysis of networks and, more generally, complex interacting systems have emerged which leverage topological techniques (\cite{Horak_2009,ShapeCollab2017,lee2011discriminative, Rieck2018}).  
These techniques generally are referred to as Topological Data Analysis (TDA) (\cite{GhristAppliedTopology,patania2017topological}). 
TDA is a relatively modern subject based on classical Algebraic Topology (\cite{Hatcher2002AlgTop,Munkres84ElAlgTop}) and that was sparked from a handful of seminal works in the late 90's (\cite{frosini1990distance,Delfinado95,Edel02Simplification, Zomorodian05Computing,Cohen-Steiner2007}). 
TDA typically endows a large variety of datasets with a notion of shape (more properly, with a topological structure) and, based on that, studies the considered data in terms of its topological features.

This field is undergoing a rapid expansion thanks to its rooting in the powerful languages of homological algebra and category theory, which provide strong formal foundations, as well as to the wide variety of applications it found, that span material science (\cite{HiraokaNakamuraHirata2016, lee2017quantifying}), biology and chemistry (\cite{chan2013topology,Meng2020WeightedPH,chung2009persistence, dequeant2008comparison, wang2005coarse, martin2010topology,phinyomark2017navigating}), sensor networks (\cite{deSilva2007coverage}), cosmology (\cite{CosmicWebPersistence}), medicine and neuroscience (\cite{patania2019topological,lawson2019persistent, GiustiPastalkovaCurto2015, WangOmbaoChung2017, YooKimAhn2016,PetriVaccarino14,ibanez2019topology,lord2016insights,ibanez2019spectral}), manufacturing and engineering (\cite{GuoBanerjee2016,phinyomark2018analysis,campbell2019differences}), social sciences (\cite{patania2017shape,benson2018simplicial}), and network science itself (\cite{BianconiTopPercol, Bianconi19Explosive,Kannan2019PersistentHO,Scolamiero13Strata,Patania17Top,donato2016persistent,SizBassClass17}). \\
The most central tool in TDA is undoubtedly Persistent Homology (\cite{GhristAppliedTopology,Edels2010Comp}). The theory of (or around) persistence has recently been proposed as a framework for the topological skeletonization of spaces, particularly weighted graphs and networks (\cite{Kurlin1D2015,kalisnik2019higher,Wang11ReebSkelet,chazal2015ReebGromov}). \\
In \cite{PetriVaccarino14}, the generators of persistent homology are used to build one instance of network skeletonization called \emph{homological scaffold}. However, the method has a serious drawback, consisting in the large degree of arbitrariness in the choice of one representative cycle from the many equivalent generating cycles of the same homology class. 
This is unfortunately a direct consequence of the homology classes being equivalence classes and affects all attempts to localize cycles (\cite{sizemore2018cliques,lord2016insights}).
In this work, we set out to address this issue by searching for a form of \emph{canonicity} in the choice of generators, namely by computing \emph{minimal representatives} of homology bases. \\
Minimal homology bases have long been investigated (\cite{Obayashi2018VolOptimal, Dey09}), with a breakthrough only coming thanks to the introduction of a first efficient algorithm for the computation of bases in dimension one (\cite{Dey18Eff}). 
Here, we leverage said minimal bases to propose a new approach to network skeletonization, the \emph{minimal scaffold}, which overcomes the limitation of the previous one. While the minimal scaffold is not unique in the most general case possible, we provide strong guarantees and caveats on when and to what degree it is well-defined.
We then show a few applications of the novel method, concluding the paper with a comparison between our and the previous construction. \\

\noindent
The paper is organized as follows. 
Section \ref{sec:background} provides a brief overview of the main concepts in Topological Data Analysis. 
Section \ref{sec:old-scaffold} describes the original approach to network skeletonization by means of persistent homology, and highlights the deficiencies which we wish to address. 
In Section \ref{sec:minimal-basis}, the topic of computing minimal representatives of a homology basis is worked out. 
Section \ref{sec:minimal-scaffold} introduces the main concept of this work, the minimal scaffold. 
In Section \ref{sec:uniqueness}, the issue of uniqueness is discussed, with some results stated, leading to a more refined version of the minimal scaffold. 
Section \ref{sec:applications} showcases some applications for the minimal scaffold. 
In the light of its computational complexity, we further carry out in Section \ref{sec:comparison} a statistical comparison between the minimal and original scaffolds, providing some heuristic guarantees and caveats. 
Section \ref{sec:conclusions} concludes the discussion. \\

\noindent
List of symbols and their common usage throughout the paper: \\

\begin{tabular}{|c|l|}
    \hline
    \textbf{Symbol} & \textbf{Meaning} \\
    \hline
    $C$ & A point cloud in $\mathbb{R}^d$ \\
    $K$ & A simplicial complex \\
    $\mathcal{F}$ & A filtration of simplicial complexes $(K^\varepsilon)_{\varepsilon=1, ..M} $\\
    $W$ & A non-negatively weighted finite graph \\
    $V$ & The set of vertices of a graph  \\
    $E$ & The set of edges of a graph \\
    $\text{VR}(W)$ & The Vietoris-Rips complex of graph $W$ \\
    $C_k(K)$ & The vector space over $\mathbb{Z}_2$ of chains of $k$-simplices of the complex $K$\\
    $\partial_k$ & The boundary operator between $C_k(K)$ and $C_{k-1}(K)$ \\
    $H_1(K)$ & The $1^{st}$ homology group of complex $K$ \\
    $\beta_1(K)$ & The dimension of $H_1(K)$ \\
    $PH_1(\mathcal{F})$ & The 1-dimensional persistent homology of filtration $\mathcal{F}$ \\
    $\mu$ & A function assigning non-negative weights to edges and cycles \\ 
    $B$ & A minimal homology cycle basis \\
    $\tilde{B}$ & A minimal homology cycle basis with draws \\
    $B^*$ & The disjoint union of minimal cycle bases across a filtration \\
    $\tilde{B}^*$ & The disjoint union of minimal cycle bases with draws across a filtration \\
    $V_i$ & A set of homologous, equally minimal variants of a basis cycle \\
    $\mathcal{H}(W)$ & The homological scaffold of weighted graph $W$ \\
    $\mathcal{H}_{min}(W)$ & The minimal homological scaffold of weighted graph $W$ \\
    $\tilde{\mathcal{H}}_{min}(W)$ & The minimal homological scaffold with draws of weighted graph $W$ \\
    \hline
\end{tabular} \\

\section{Background}
\label{sec:background}

In this section we introduce the minimum amount of mathematics necessary to the understanding of the rest of the paper. We refer to classical textbooks on the subject for further reading (\cite{Hatcher2002AlgTop,Munkres84ElAlgTop,Edels2010Comp,GhristAppliedTopology}).
\subsection*{Simplicial complexes}
Thanks to their proven flexibility in a plethora of applicative contexts, simplicial complexes are the most adopted mathematical structure for encoding unorganized, large-size and high-dimensional data.
In purely combinatorial terms, a (finite) {\em simplicial complex} $K$ on a finite set $V$ is a collection of non-empty subsets of $V$, called {\em simplices}, with the property of being closed under inclusion, i.e., every non-empty subset of a simplex of $K$ is itself a simplex of $K$.
Given a simplicial complex $K$, the elements of $V$ are called {\em vertices} of $K$ and a simplex $\sigma \in K$ is called a {\em $k$-simplex} (equivalently, a simplex of dimension $k$) if it consists of $k+1$ vertices.
The {\em dimension} of a simplicial complex $K$ is the largest dimension of the simplices in $K$.
\\
Even if the abstract definition of a simplicial complex just given is able to capture a variety of datasets not necessarily endowed with a geometrical realization, it is worth to be mentioned that, intuitively, a simplicial complex is nothing but a collection of well-glued bricks, its simplices. According with such a perspective, a $k$-simplex can be seen as the convex hull of $k+1$ geometrically independent points.
For instance, a 1-simplex is an edge, a 2-simplex is a triangle, a 3-simplex is a tetrahedron, and so on.\\

\subsection*{Homology}
Homology is a topological tool which provides invariants for shape description and characterization.
Given a simplicial complex $K$, it is possible to associate to it a collection of vector spaces $C_k(K)$ over a field, in our case $\mathbb{Z}_2$, whose bases are indexed by the $k$-simplices so that, loosely speaking, we say that these spaces are generated by the $k$-simplices of $K$. These spaces are connected by boundary operators $\partial_k:C_k(K) \rightarrow C_{k-1}(K)$ mapping each $k$-simplex $\sigma$ in the sum of the $(k-1)$-simplices of $K$ strictly contained in $\sigma$.

We denote as $Z_k(K):=\ker\partial_k$ the space of the $k$-cycles of $K$ and as $B_k(K):=\im\partial_{k+1}$ the space of the $k$-boundaries of $K$. Then, since $\partial_{k}\partial_{k+1}=0$, the quotient $$H_k(K):=\frac{Z_k(K)}{B_k(K)}$$ defines a vector space called \textit{$k^{th}$ homology group} of $K$.

We will call two $k$-cycles {\em homologous} if they belong to the same homology class.
\\
Roughly speaking, homology reveals the presence of ``holes" in a shape.
A non-null element of $H_k(K)$ is an equivalence class of cycles that are not the boundary of any collection of $(k+1)$-simplices of $K$.
Such classes represent, in dimension 0, the connected components of complex $K$, in dimension 1, its tunnels and its loops, in dimension 2, the shells surrounding voids or cavities, and so on.\\

\subsection*{Persistent homology}
An intrinstic limitation of homology concerns the need for working with a single simplicial complex representing the dataset under investigation. 
However, in real world applications, the presence of noise and of measurement errors makes the choice and construction of a single steady representation very hard in practice.  
\textit{Persistent homology} (\cite{Edel02Simplification,Edels2010Comp}), currently one of the main tools in Topological Data Analysis, aims at solving this issue through a multi-scale study of a dataset and of its homological features by associating to it a sequence of simplicial complexes.
The concept of filtration captures exactly the idea of analyzing a dataset at different thresholds of a parameter on which it depends. 
More formally, given a simplicial complex $K$, a {\em filtration} $\mathcal{F}$ of $K$ is a sequence of its subcomplexes such that
$$\emptyset\subseteq K^1 \subseteq \dots \subseteq K^M=K$$
Given a filtration of a simplicial complex $K$, persistent homology keeps track of the evolution of the non-null non-homologous cycles of $K$ and, associating a lifespan to each of them, is able to discriminate the relevant information from the noise.
Formally, for $p,q=1, \dots, M$ with $p<q$, $H_k^{p,q}(\mathcal{F})$ on $(p,q)$ of a filtration $\mathcal{F}$ consists of the image of the linear map between $H_k(K^p)$ and $H_k(K^{q})$ induced by the inclusion of complexes between $K^p$ and $K^q$. So, more intuitively, the elements in $H_k^{p,q}(\mathcal{F})$ represent the cycles of $K$ which survive from step $p$ to step $q$.\\
Given a filtration of finite simplicial complexes $\mathcal{F}$, we define its $k$-dimensional {\em persistent homology classes} as the homology classes of $\bigoplus_\varepsilon H_k(K^\varepsilon)$ modulo the maps induced by the inclusion of simplicial complexes. More properly, $h_1\in H_k(K^p)$ and $h_2\in H_k(K^q)$ with $p\leq q$ are equivalent if and only if $\iota^{* \ p,q}_k(h_1)=h_2$ where $\iota^{* \ p,q}_k$ denotes the linear map between $H_k(K^p)$ and $H_k(K^{q})$ induced by the inclusion of complexes between $K^p$ and $K^q$. We call \emph{$k$-dimensional persistent homology} $PH_k(\mathcal{F})$ the space spanned by the $k$-dimensional persistent homology classes. \\
As proven in \cite{Zomorodian05Computing}, a basis of $PH_k(\mathcal{F})$ is in bijective correspondence with a finite set of intervals of the form $\{(p,q), \ p<q, \ p,q \in \mathbb{Z}\cup\infty \}$ referred as {\em persistence pairs}. We define a set of $k$-dimensional {\em generator cycles} of the persistent homology as a set of $k$-cycles of $K^M$ whose persistent homology classes form a basis of $PH_k(\mathcal{F})$. \\
The information about the ``life" of each homology class can be collected in a visual, informative representation of the topological structure of the input, the {\em persistence barcode}: a plot consisting of a bar for each homological feature appearing throughout the filtration, stretching from its birth to its death value.
An equivalent way to depict the same information is through the \emph{persistence diagram}: the persistence diagram is the multi-set (i.e., multiple instances of the same element are allowed) of points in $\mathbb{R}^2$ consisting of all the \emph{(birth, death)} pairs, i.e., pairs of values $p < q$ such that a $k$-dimensional homology class arises at filtration step $p$ and becomes zero at step $q$.
\noindent
Persistent homology owes its popularity as a descriptor to the immediacy and power of these visual representations of the homological information but, even more, to the fact that the retrieved features are provably stable. 
In fact, by defining a notion of distance among persistence diagrams or barcodes, it can be shown that similar datasets necessarily have similar homological features (\cite{Cohen-Steiner2007}).

\subsection*{Building (filtered) complexes}
\begin{figure}
    \centering
    \includegraphics[scale=0.43]{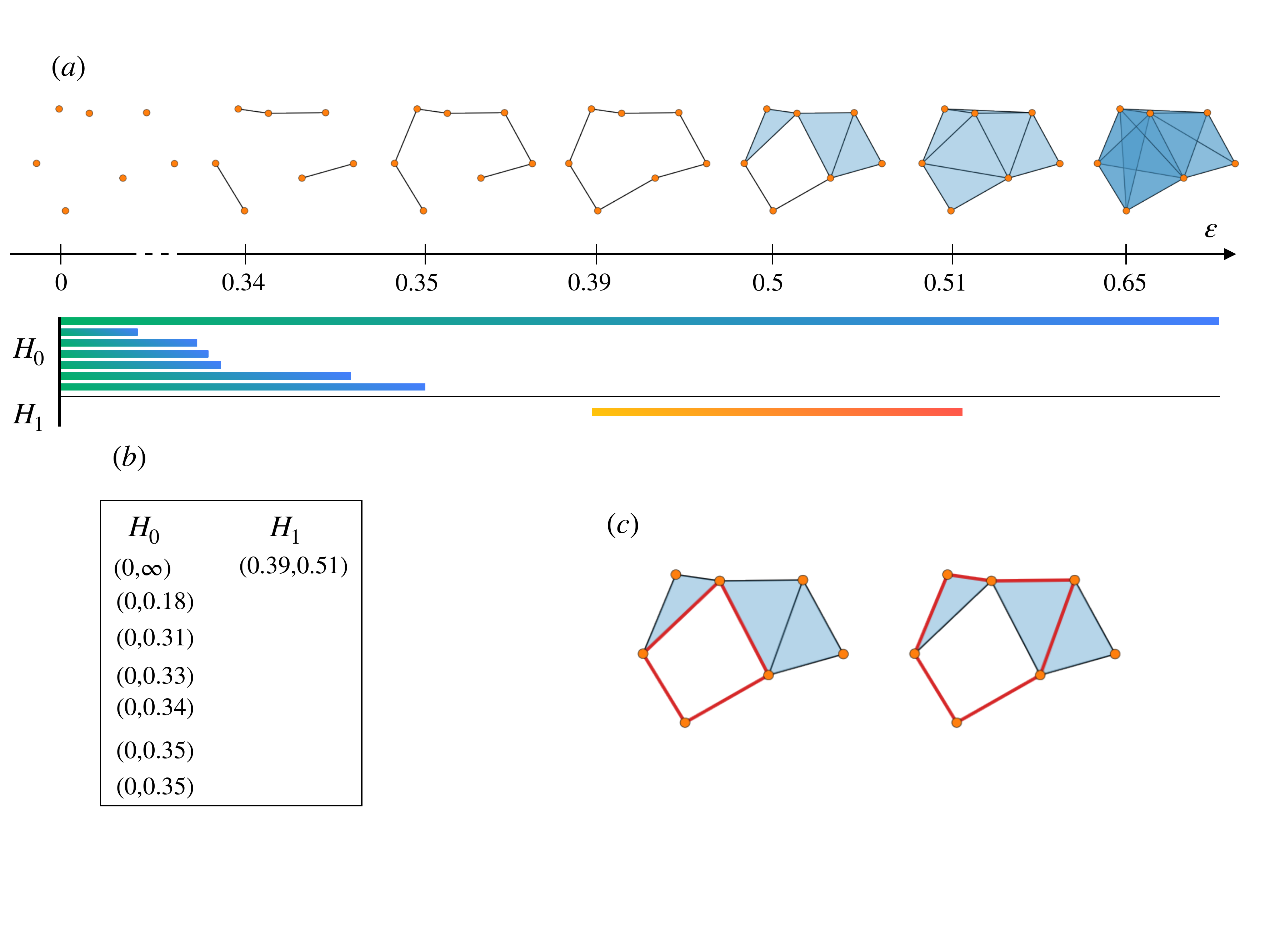}
    \caption{(a) An example of Vietoris-Rips filtration of simplicial complexes with parameter $\varepsilon$, and the corresponding barcode for 0- and 1-dimensional persistent homology. (b) The persistent pairs of the above filtration. (c) Two equivalent representatives of the (only) generator of $PH_1$.}
    \label{fig:FiltrationPanel}
\end{figure}

In many applications, one is not directly called to deal with a simplicial complex, but has instead access to data in the form of point clouds in a metric space or of  weighted graphs. 
For example, data may be obtained as a sample of some (unknown) ground truth, i.e., an undisclosed manifold of dimension usually much lower than the space it is embedded in (\cite{GhristAppliedTopology}).
Another typical subject of application is network science (\cite{Scolamiero13Strata,SizBassClass17}): 
in this setting, the input is in the form of a weighted graph. 
Notice that in this case it is not mandatory that the graph can be embedded in some metric space, i.e., that the edge weighting respects a triangular inequality. 
Networks are not necessarily representations of geometrical entities, and still the topological approach extends naturally to this context.  \\
In both these cases, one needs to provide a suitable simplicial complex resting on the given structure. The subject has been addressed extensively (see, for example, \cite{Edels2010Comp}); in here, we simply review the most typical scheme, called the \textit{Vietoris-Rips complex}. Given a graph $G = (V,E)$, its \emph{flag} or \emph{clique} complex is the simplicial complex $Flag(G)$ whose simplices coincide with the cliques of $G$.
\\
Given a point cloud $V \subset \mathbb{R}^n$ and fixed a value $\varepsilon > 0$, one can build a graph $G^\varepsilon$ with a vertex for every point in $V$, and an edge between two vertices every time the distance between the corresponding points is less or equal than $\varepsilon$.
Analogously, given a weighted graph $G=(V,E)$ one can build a subgraph $G^\varepsilon$ on the same vertex set, with only those edges that have weight less or equal than $\varepsilon$.
Independently from the considered case, one can define the {\em Vietoris-Rips complex} $VR^\varepsilon$ of parameter $\varepsilon$ as the flag complex $Flag(G^\varepsilon)$ of graph $G^\varepsilon$. Furthermore, since varying $\varepsilon$ the Vietoris-Rips complexes $VR^\varepsilon$ form an increasing sequence of simplicial complexes, the family $(VR^\varepsilon)$ gives raise to a filtration denoted as {\em filtered Vietoris-Rips complex} (see Fig. \ref{fig:FiltrationPanel}).
\\
As already mentioned, Vietoris-Rips complexes are employed in a wide variety of different application domains. 
The reason is that their definition only depends on the pairwise distances between points, making them efficient to compute and to store with respect to more refined alternatives. It is worth noticing, however, that cost of this simplicity is the fact that the dimension of a Vietoris-Rips complex can explode even when the points are sampled from a low-dimensional subspace of $\mathbb{R}^n$.

\FloatBarrier

\section{Homological Scaffold}
\label{sec:old-scaffold}
The homological scaffold originated from the intuition that traditional, graph-theoretical tools in network analysis were naturally able to capture significant properties (\cite{baronchelli2013networks}), but proved not as effective in detecting multi-agent and large-scale interactions. 
Interest in searching for alternative descriptors of network relations arose, and soon works were published which leveraged invariants offered by computational topology (\cite{lum2013extracting, lee2011discriminative, ShapeCollab2017}). \\
In proposing the scaffold (\cite{PetriVaccarino14}), the authors pointed out that homological might be able to summarize well network \emph{mesoscale} structures, i.e., features living between the purely local connections and the global statistics, to which previous methodologies were blind. 
Furthermore, this structure could be analyzed over the continuous, full range of interaction intensities, without the need for ad-hoc domain-specific thresholds. \\
Homological cycles intuitively describe obstruction patterns. 
The presence of non-trivial homology within a given region of a network highlights its structure as non-contractible, binding signals to flow over constrained channels, which in turn play the role of bridges.  \\
To test the method, the homological scaffold was computed from resting-state fMRI data for 15 healthy volunteers who were either infused with placebo or psilocybin: the scaffold discriminated the two groups, as well as providing meaningful insight as to the impact of the psychoactive substance onto the pattern of information flow in the brain \cite{PetriVaccarino14}. \\

\noindent
Given a non-negatively weighted finite graph $W = (V,E,w \ : \ E \mapsto \mathbb{R}^+$), let $\mathcal{F}$ be a filtration of simplicial complexes as above. \\
Let $\{ b_i \}$ be a set of $1$-dimensional generator cycles of the persistent homology. Since we are over $\mathbb{Z}_2$, each of the $ b_i$'s is completely identified by its support, which is a set of edges of $E$. In particular, we can depict set $\{ b_i \}$ as a matrix whose rows are indexed by $E$ and having the $b_i$'s as columns. The row sums, as natural numbers, form a new weighting function on the edges of $W$, the new weights counting precisely in how many persistent cycles an edge appears along the filtration.

\begin{definition} 
Suppose $W$ and $\mathcal{F}$ as above, and consider a set $\{ b_i \}$ of $1$-dimensional generator cycles of the persistent homology. Consider the function $h_W : E \mapsto \mathbb{R}^+$

\begin{equation} \label{defHomScaff}
    h_W := \sum_{i} \mathbbm{1}_{e \in b_i}
\end{equation}
where by $\mathbbm{1}_{e\in b_i}$ we denote the indicator function $E \mapsto \mathbb{R}^+$ such that $\mathbbm{1}_{e\in b_i}(e')=1$ if $e'$ appears in $b_i$, and 0 otherwise.

\noindent
Then the \textbf{homological scaffold} of $W$ is the weighted graph $\mathcal{H}(W)$ such that
\begin{itemize}
\item[-] its vertex set coincides with the vertex set of $W$
\item[-] its edge set $E_{\mathcal{H}}$ is a subset of the edge set of $W$, consisting of edges with nonzero value for $h_W$
\item[-] its weight function is the restriction of $h_W$ to $E_{\mathcal{H}}$. 
\end{itemize}
\end{definition}

\noindent
In accordance with the above definition, building the homological scaffold of a weighted network $W$ is a method of \emph{network compression} or \emph{skeletonization}. The definition also implies that edge weights are assigned by the number of basis cycles the edge belongs to. \\

\begin{figure}
    \centering
    \includegraphics[scale=0.43]{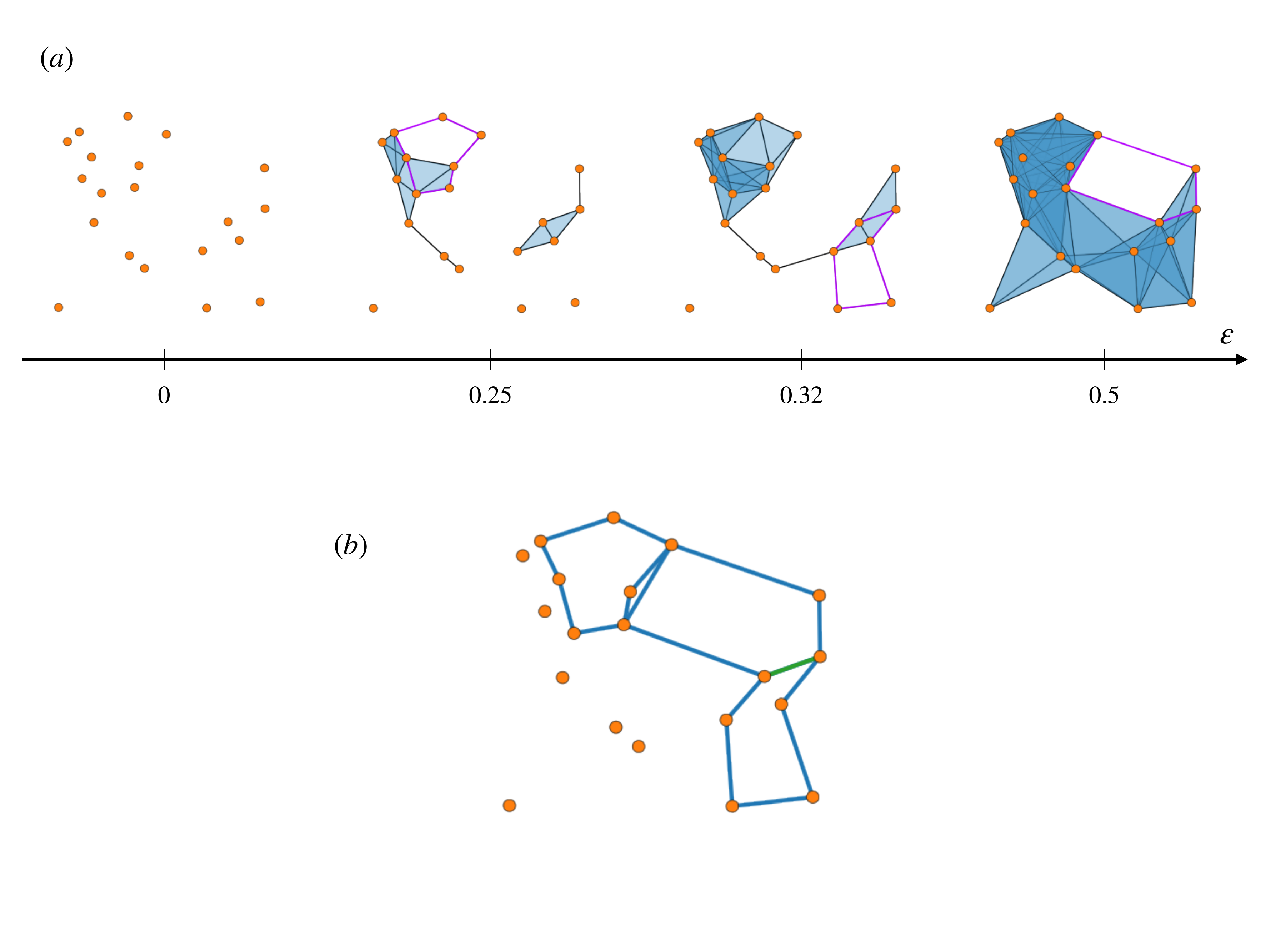}
    \caption{(a) A point cloud in $[0,1]^2$ and the generators of $PH_1$, plotted on the filtration step they appear at (scale reported on the axis below). (b) The resulting homological scaffold. Edges in blue have weight 1, each belonging to only one generator. The edge in green has weight 2, as it belongs to two generators.}
    \label{fig:HomScaffold}
\end{figure}

\noindent
We provide an example, referring to Fig. \ref{fig:HomScaffold}. In panel (a), a filtration of simplicial complexes arising from a point cloud is depicted. At each step, highlighted in purple is a representative of a persistent cycle (i.e. of a bar in the barcode), each at the scale at which it is born. \\
In panel (b), the corresponding homological scaffold is represented: it amounts to taking the union of the cycles of panel (a), i.e. stacking generators of $PH_1$, each contributing unitary weight. \\

\noindent
In the following, we shall sometimes refer to the homological scaffold as the \emph{loose}, or \emph{original} scaffold, to contrast it with the new definition of scaffold to follow. \\

\FloatBarrier

\begin{figure}
    \centering
    \includegraphics[scale=0.43, trim=0 200 0 30,clip]{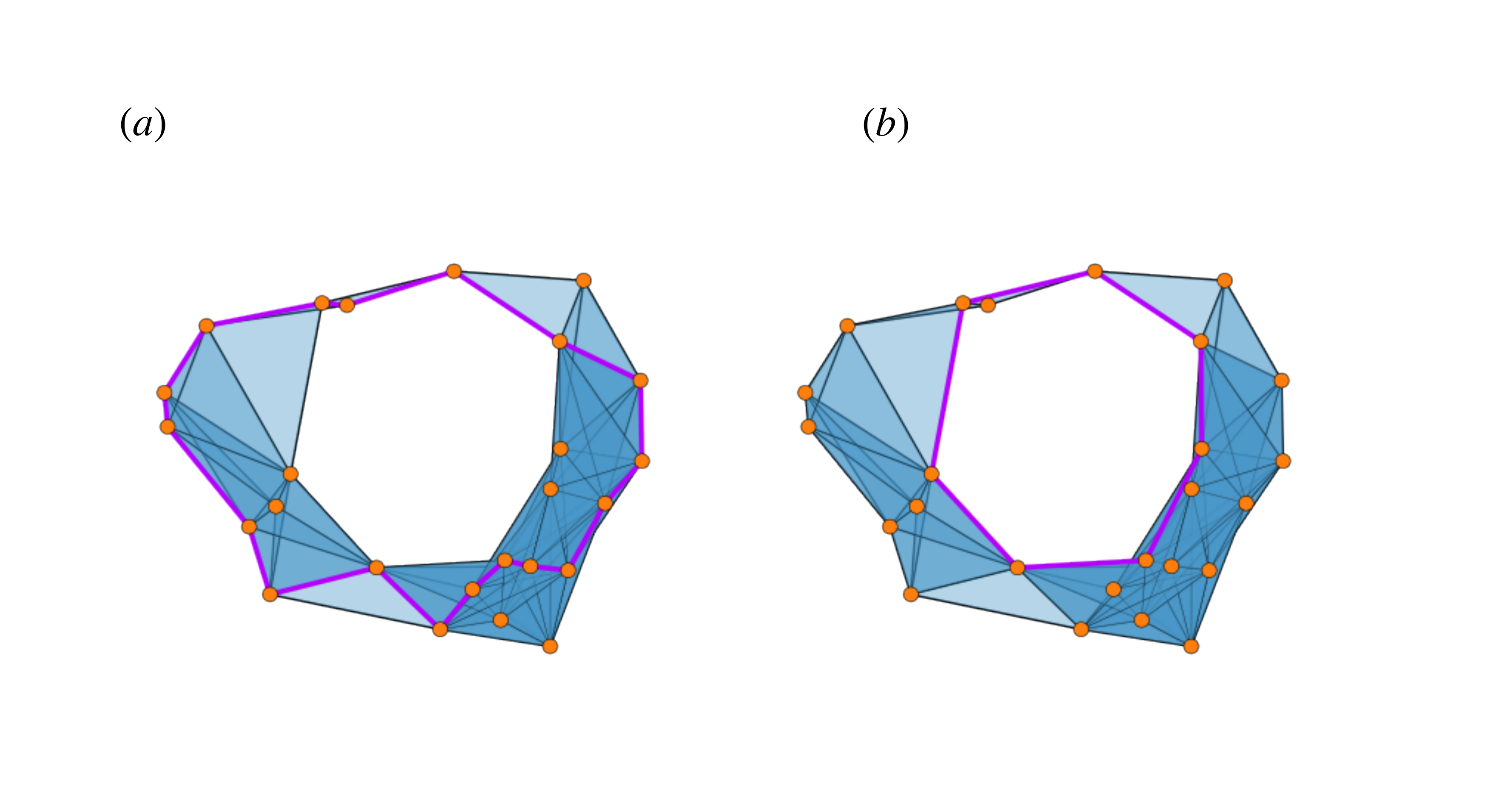}
    \caption{A simplicial complex $K$ with $\dim H_1(K) = 1$. Its homological scaffold (on a subset of the filtration steps, for clarity) is reported in panel (a): the chosen generator meanders around the hole. Furthermore, a different ordering of the list of simplices fed to the algorithm could return a different cycle. In panel (b), the shortest representative cycle is chosen: this choice is stable with respect to any ordering of the input, while at the same time endowing the generator with some metric and geometric meaning. }
    \label{fig:Cycles}
\end{figure}

\noindent
As anticipated in the introduction, it is apparent that there is a substantial source of arbitrariness in this definition. \\
Several different representative cycles exist which form a basis of the persistent homology (as a consequence of several different cycles belonging to the same homology class), and hence one must make a choice. For example, Fig. \ref{fig:Cycles}(a) depicts one specific cycle whose homology class generates (part of) the persistent homology group of the point cloud. At the same time, any other choice of edges forming a cycle around the hole is homologically equivalent and, in principle, legitimate. \\
In the original paper, the authors resorted to using the cycles as output by the \emph{JavaPlex} implementation (\cite{Javaplex}) of the persistent homology algorithm (based on the original implementation of \cite{Delfinado95}), and a posteriori checked the selected cycles for consistency. 
However, in principle, this means that the same simplicial complex written with two different orderings of the simplices could lead to different choices of generators, and therefore, to different scaffolds. \\
As such, we must be careful in the choice of nodes and edges output by the algorithm; while the presence of a generator denotes undeniably that an obstruction pattern exists, we cannot be as confident about its precise location in the network or the constituents that provide bridges around it. 
The homological scaffold defined in this way introduces noise in the localization of mesoscale patterns onto individual nodes and edges, a process which, if accurate, could provide valuable insight as to the functional role of single players in a network.  \\
In this work, we try to work around the problem of cycle choice and give a stricter definition, by requiring that, among all possible representatives, those of \emph{minimal total length} are chosen (e.g., Fig. \ref{fig:Cycles}(b)). \\
The original algorithm reported a computational complexity of the order $O(n^3)$ to obtain representatives of basis cycles. \\

\FloatBarrier

\section{Minimal Bases}
\label{sec:minimal-basis}

The search for minimality in the computation of the scaffold was made feasible by the introduction of efficient algorithms to compute the minimal representatives of a homology bases in dimension one. \\
It is known that in dimension higher than one, minimal representatives of a homology basis will remain elusive. 
Indeed, Chen and Freedman (\cite{ChenFreedman2011}) proved that the problem of obtaining these minimal representatives is computationally intractable, being at least as hard as the notoriously NP-Hard Nearest Codeword Problem. 
Furthermore, it is even NP-Hard to approximate within any constant factor, meaning that no polynomial-time algorithm exists to obtain an approximate minimal basis that differs from the exact one by at most a multiplicative constant. In the light of this, we must necessarily restrict our attention to the 1-dimensional case, i.e., computing minimal representatives of a basis of $H_1$.\\

\subsection{Minimal Bases and Dey's Algorithm}

Given a simplicial complex $K$, let us consider $C_1$ the vector space generated by the 1-simplices of $K$ and $Z_1$ the vector space of $1$-cycles, i.e., $Z_1= \ker \partial_1$.
Given a $1$-cycle $b \in Z_1$, let $\mu(b)$ be its length, i.e., the sum of the weights of the $1$-simplices that form it, and denote by $[b]$ the homology class $b$ belongs to. Finally, let $\beta_1 := \dim H_1(K)$. We want to obtain a set of $\beta_1$ $1$-cycles $\in Z_1$

\begin{equation} \label{eqMinBasis}
    \{ b_1 , ... , b_{\beta_1} \} = \underset{\text{Span}\{ [b_i] \} = H_1 }{\text{argmin}} \sum_i \mu \left( b_i \right)
\end{equation}

that is a set of cycles of minimal length whose homology classes span $H_1(K)$. In accordance with the literature, we call this set a \emph{minimal homology basis},  with a slight abuse of terminology, as it would be more appropriate to call it a \emph{minimally-represented homology basis}. \\
In 2018, Dey et al. (\cite{Dey18Eff}) introduced a polynomial-time algorithm to obtain said representatives. 
Building on the work of Horton (\cite{Horton87}), de Pina (\cite{dePina95}), and Mehlhorn et al. (\cite{Mehlhorn04}), the algorithm sets off to compute a basis of the space of cycles. 
Then, it applies a cohomological technique called \emph{simplex annotation} (\cite{Dey12Annotating}) to lift a basis of cycles to a basis of the homology group $H_1$, while at the same time enforcing the minimal length constraint.
A sketch of the algorithm follows. \\
\textsc{Algorithm: MinBasis($K$)}
\begin{itemize}
    \item A basis of the cycles group $Z_1$ is found via a spanning tree. Each edge in the complement of the spanning tree identifies a candidate cycle (\cite{Horton87}).
    \item An annotation of the edges is computed via matrix reduction (\cite{Dey12Annotating}). This yields the dimension $\beta_1$ of $H_1$, as well as an efficient tool to determine if two cycles $b_1$ and $b_2$ are linearly dependent in $H_1$ ( $[b_1] = [b_2]$).
    \item A set of \emph{support vectors} is generated which maintains a basis of the orthogonal complement in $H_1$ of the minimal basis cycles.
    \item Iteratively for each dimension of $H_1$, the candidate set of cycles is parsed in search of cycles $b$'s that are linearly independent in homology from the previous ones (exploiting the support vectors). Among these, the $\mu$-shortest one is added to the minimal basis.
    \item The set of support vectors is updated for the remaining dimensions to enforce it remain a basis of the orthogonal complement of the basis. \item The last two steps above are repeated until completion of the minimal basis. 
\end{itemize}
\vspace{10pt}
Call $B = \{ b_i \}$ the output of \textsc{MinBasis} on input $K$. \\
\textbf{Theorem} (3.1, \cite{Dey18Eff}) Cycles in $B$ form a minimal homology basis of $H_1(K)$. \\

\noindent
Notice that the minimal homology basis is guaranteed to exist, as we only work with finite simplicial complexes, which imply the existence of a finite number of bases. 
However, it needs not, in general, be unique. Several different cycles of the same minimal length may all belong to the same homology class of a basis cycle. 
Heuristically, this is especially true in case the input complex is unweighted (equivalently, has equal weights for every edge), in which case the length of a cycle is the number of edges that form it. 
Furthermore, there exist cases when different sets of cycles of minimal length generate the same homology space, and are not even pairwise homologous. 
We will treat the problem of the uniqueness of the minimal basis in more detail in the following, and account for it explicitly in the construction of the minimal scaffold. \\

\noindent
The computational complexity of the above procedure is evaluated (\cite{Dey18Eff}) to $O(n^2 \beta_1 + n^\omega)$ where $n$ is the number of simplices in $K$ and $\omega$ is the fast matrix multiplication exponent, which as of 2014 is bounded by 2.37 (\cite{Dey18Eff, Coppersmith1990, LeGall2014}). This yields a worst-case complexity of $O(n^3)$ in the number of simplices for general complexes, which we recall is itself of order $3$ in the number of points in the worst case. \\

\section{Minimal Scaffold}
\label{sec:minimal-scaffold}
\FloatBarrier
In this section, we introduce an alternative definition for the homological scaffold, which we call minimal, based on the minimal representatives obtained above, and aims at overcoming the arbitrariness in the cycle choice of the previous definition. 
After addressing the simplest case, we analyze its uniqueness properties and introduce a second, more refined, definition.

Let $\mathcal{F}$ be the filtration of simplicial complexes induced by a non-negatively weighted finite graph $W$. For all filtration steps $\varepsilon$, define, as per (\ref{eqMinBasis}), $B^\varepsilon := \{ b_i^\varepsilon \}$ the minimal homology basis of $H_1(K^\varepsilon)$. Take the disjoint union of minimal bases for $\varepsilon$ varying on all filtration steps
\[ B^* := \coprod_\varepsilon B^\varepsilon \]

\begin{definition}
Suppose $W$, $\mathcal{F}$ and $B^*$ as above. Similarly to the loose case, define the function $h_{W, min} : E \mapsto \mathbb{R}^+$ as
\begin{equation} \label{defMinScaf}
    h_{W, min} := \sum_{b  \in B^*} \mathbbm{1}_{e \in b}
\end{equation}
Then, we define the \textbf{minimal scaffold} of $W$ as the weighted graph $\mathcal{H}_{min}(W)$ whose:
\begin{itemize}
\item[-] vertex set coincides with the vertex set of $W$
\item[-] edge set $E_m$ is a subset of the edge set of $W$, consisting of edges with nonzero value for $h_{W, min}$
\item[-] weight function is the restriction of $h_{W, min}$ to $E_m$. 
\end{itemize}
\end{definition}

\noindent
The minimal scaffold amounts, again, to the stacking of generator cycles across a filtration. However, two differences are to be noted with respect to the loose definition.
First, we require the representative cycles to be minimal. 
Second, we point out that while the loose scaffold is built by aggregating the generator cycles of $PH_1(\mathcal{F})$, the minimal scaffold is built by independently computing a minimal basis for each $H_1(K^\varepsilon),$ for all $\varepsilon$.
Notice that, since cycles are modified throughout a filtration, it would be meaningless to talk about a minimal representative over a certain persistence interval. 
This also means that its computation can be effectively parallelized by assigning different filtration steps to different jobs, and later recombining the outputs. \\
An interesting phenomenon that descends directly from the above peculiarity is that the minimal scaffold of random point clouds tends to display a more pronounced  triangular structure (clustering) around cycles. 
Indeed, as longer (or, in non-metrical filtrations, later) edges are introduced, a cycle can be shortened (by the triangular inequality) by a longer edge which cuts a corner. 
Since at each step the algorithm records the minimal representative, upon aggregating the minimal scaffold one finds each cycle in its progressively shorter version, and the \emph{history} of the shortening is visible as a padding of triangles around it.  \\

\begin{figure}
    \centering
    \includegraphics[scale=0.43
    ]{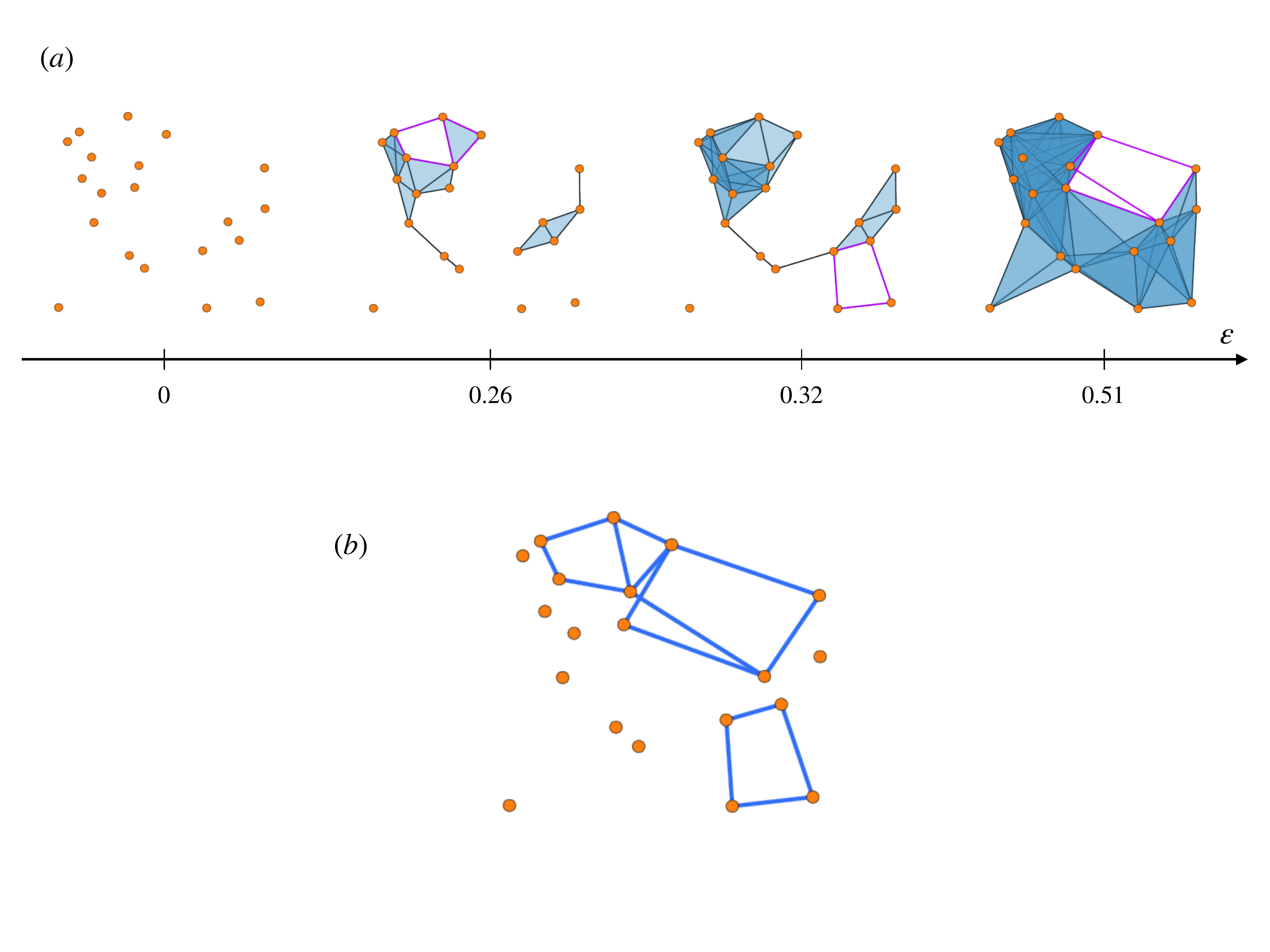}
    \caption{(a) The same point cloud of Fig. \ref{fig:HomScaffold}. Along the filtration we show the evolution of minimal generators, which can get progressively shorter as new edges are introduced. For example, at $\varepsilon = 0.26$, the pentagonal cycle gets cut to a shorter quadrilateral, albeit with an individual longer edge. This evolution is accounted for in the minimal scaffold, which displays the triangle-rich structure mentioned above.  (b) The resulting minimal scaffold (weights not reported).}
    \label{fig:MinScaffoldExample}
\end{figure}

\noindent
Considering the example of Fig. \ref{fig:MinScaffoldExample}, in panel (a) we observe an example of a filtration of simplicial complexes. At each step, highlighted in purple we may see the minimal representative of a homology class, together with its evolution history. At filtration value $0.26$, we observe a pentagon being reduced to a shorter, quadrilateral cycle by the addition of a longer edge. This is an example of the phenomenon explained above. Fig. \ref{fig:Cycles} gives a visual description of the difference between a minimal and generic cycle. \\
The union of these progressively shorter cycles for all steps (weighted according to Definition \ref{defMinScaf}) is the minimal scaffold, as seen in Fig. \ref{fig:MinScaffoldExample} panel (b). \\

\noindent
We remark that, if there is no ambiguity in the construction of a filtration of simplicial complexes from a point cloud, or from a weighted graph, we will indifferently speak of the scaffold as a function of either of them ($\mathcal{H}_{min}(C)$, or $\mathcal{H}_{min}(W)$, or $\mathcal{H}_{min}(\mathcal{F})$). \\

\noindent
We have mentioned that the scaffold amounts to a change in weighting in the input graph 
\[ h_{W,min} \ : \ E \mapsto \mathbb{R}^+ \]
altering the original weights of the edges. Additionally, considering \emph{node strength} (i.e. the sum of the weights of the edges incident to a given node), it can equally be considered as a function
\[ \mathcal{H}_{min} \ : \ V \mapsto \mathbb{R}^+ \]
assigning weights to nodes. 
Considering the reliability of the choice of edges in the procedure, this explains why the minimal scaffold can be utilized to associate mesoscopic features with single nodes and links. \\

\FloatBarrier
\subsection*{Computational Complexity}
For large input sizes, the cost of assembling the minimal basis cycles into the scaffold is negligible with respect to the cost of computing such minimal basis. 

\begin{figure}
        \centering
        \includegraphics[scale = 0.6]{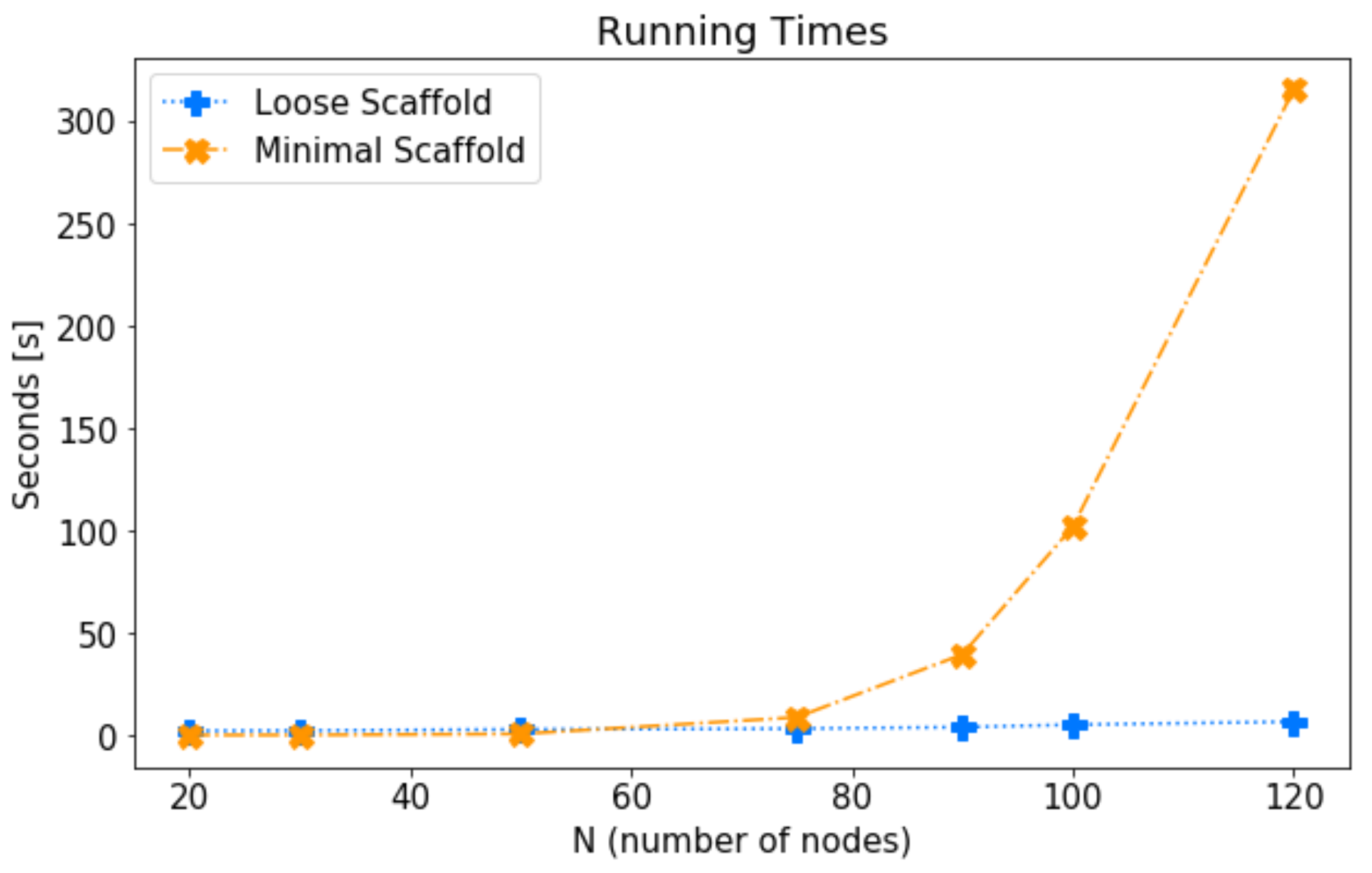}
        \caption{The running times of computing the minimal and loose scaffolds for Watts-Strogatz weighted random graphs. For all instances, number of nodes $N$ is indicated on the x-axis. Number of stubs $k$ is $N/2$, and rewiring probability is $p=0.025$.}
        \label{fig:RunningTimes}
    \end{figure}

We know that each run of Dey's algorithm costs $O( \vert K \vert^3)$ in the worst case (\cite{Dey18Eff}), and in the worst case $\vert K \vert$ is itself $O(n^3)$ where $n$ is the number of points. \\
The number of filtration steps has an upper bound of $O(n^2)$ (i.e., the number of edges) in the worst case, as in general every edge may carry a different weight. Hence Dey's algorithm has to be run once for each edge in the worst case. \\
This yields a theoretical worst-case complexity of order $O(n^9 n^2) = O(n^{11})$. Therefore, while the minimal scaffold is undeniably a polynomial-time algorithm, its practical computation is often hindered by its dire lack of scalability, especially if compared against the loose version, which has a far more favourable complexity. \\
A comparison of running times is carried out in Fig. \ref{fig:RunningTimes}, which clearly shows that computing the minimal scaffold on an ordinary machine can quickly become troublesome.

\FloatBarrier

\subsection*{Implementation}
\noindent
We have written a Python implementation of Dey's algorithm, together with a library for the computation of the minimal scaffold. The code is available on GitHub at \cite{GitHubMinScaf}, with some usage examples. It allows for shared-memory multi-threaded parallelism across filtration steps to improve computation times, while still being suitable for ordinary desktop workstations.

\section{Uniqueness of the minimal scaffold}
\label{sec:uniqueness}

The uniqueness of the minimal scaffold depends on the uniqueness of the minimal basis. Indeed, if there exists only one possible set $B^*$ of cycles forming a minimal basis, then the scaffold is uniquely determined.
Two issues affect the uniqueness of set $B^*$.
\subsection*{Draws} The first one arises when two or more different and homologous basis cycles are of the same minimal length. 
This case is relatively simple to work around: we modify the definition of minimal scaffold to keep track of all variants of minimal basis cycles, dividing the weight equally among them. \\
Specifically, to account for this issue we have slightly modified Dey's algorithm. In its last step described above, one is concerned with finding all cycles whose annotation is not orthogonal to the given support vector: among these, the one with minimal length is chosen as a basis cycle. Instead, we keep track of \emph{all} such cycles with the same minimal length. This does not alter the complexity, as one needs to check all possible cycles anyway. We call this case a \emph{draw}. \\

\noindent
Therefore, we modify set $B$ to become a set of \emph{sets of cycles}. 
Given complex $K$, we define a minimal basis \emph{with draws}
\[ \tilde{B} \ := \ \bigcup_{i=1}^{\beta_1(K)} \ \{ b_{i,1} , ... , b_{i, n_i} \} \]
where for all $i = 1, ..., \beta_1(K)$, the cycles $ b_{i,j}$ with $j=1, ..., n_i $ are homologous and have the same minimal length. Furthermore, for every choice of $j_i \in \{ 1 , ... , n_i\}$, $\text{Span}_i \{ b_{i,j_i} \} = H_1(K)$. Call $V_i := \{ b_{i,1} , ... , b_{i, n_i} \}$ each set of draws, i.e., \emph{variants} of the $i^{th}$ minimal basis cycle, $\forall i=1,...,\beta_1(K)$. \\
In the example of Fig. \ref{fig:PanelDraws}(a) and (b), we have set $\tilde{B} = \{ \ \{ b_{1,1} , b_{1,2}\} \ \}$, whereas set $B$ might have indifferently been equal to $\{ b_{1,1} \}$ or to $\{ b_{1,2} \}$, whichever happened to come first in the search. \\
The minimal scaffold is modified accordingly. 
Given the usual filtration $\mathcal{F}$, let $\tilde{B}^\varepsilon$ be the minimal basis with draws of $H_1(K^\varepsilon)$. 
Again, we aggregate all variants of minimal basis cycles along the filtration 
\[ \tilde{B}^* \ := \ \coprod_\varepsilon \tilde{B}^\varepsilon \]

\noindent
Then, we define the weighting function \emph{with draws} $\tilde{h}_{W, min} \ : \ E \mapsto R^+$

\begin{equation} \label{defMinScafDraws}
     \tilde{h}_{W,min} := \sum_{V \subset \tilde{B}^*} \frac{1}{\vert V \vert}\sum_{b\in V} \mathbbm{1}_{e \in b}
\end{equation}

\noindent
and the resulting {\em minimal scaffold with draws} $\tilde{\mathcal{H}}_{min}(W)$ is built from $\tilde{h}_{W,min}$ as in Definition \ref{defMinScaf}.

\noindent
The meaning of the above definition is that all variants of all minimal basis cycles are taken into account when building the scaffold, and the weights are assigned dividing each variant's contribution by its cardinality, for each filtration step. 
In the example of Fig. \ref{fig:PanelDraws}(c), the two cycles forming the variant of the only generator are multiplied by a factor of $\frac{1}{2}$ and then summed: therefore, common edges outside the diamond are assigned weight $1$, consistently with the minimal scaffold in definition (\ref{defMinScaf}), whereas the four edges forming the perimeter of the diamond each get assigned weight $\frac{1}{2}$. \\

\noindent
With the introduction of draws, we settle the case when ambiguity arises among individual cycles, without interactions. As an example, we can state the following result. \\
\textbf{Proposition} If $\mathcal{F}$ is such that, for all $\varepsilon$ in the filtration, each basis cycle belongs to a different connected component of $K^\varepsilon$, then the minimal scaffold with draws $\tilde{H}_{min}(\mathcal{F})$ is unique.

\subsection*{Pathological cases} 
\begin{figure}
    \centering
    \includegraphics[scale=0.43]{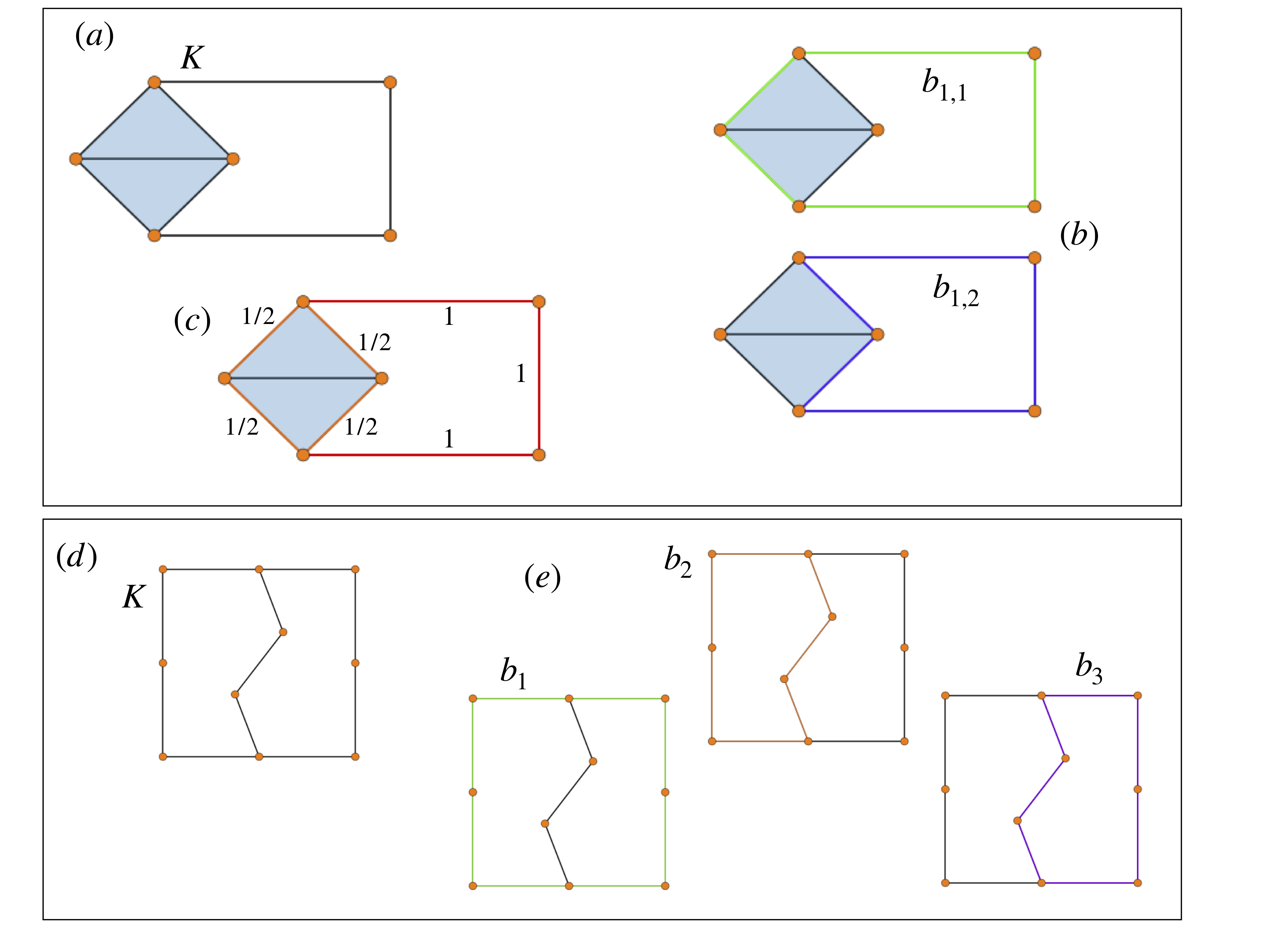}
    \caption{Top panel: (a) A simplicial complex $K$. (b) Two homologous and equally minimal generators of $H_1(K)$. (c) The minimal scaffold with draws $\tilde{\mathcal{H}}_{min}(K)$. The weight is equally divided among the variants of the minimal representative. Bottom panel: (d) A simplicial complex $K$ on the represented point cloud. $H_1(K)$ has dimension 2. (e) $\mu(b_1) < \mu(b_2) = \mu(b_3)$. A minimal basis can either be composed of $\{ b_1,b_2 \}$ or $\{ b_1,b_3 \}$, hence it is not unique.}
    \label{fig:PanelDraws}
\end{figure}

The other issue arises when there exist sets of minimal cycles that are not linearly independent. Suppose that three different cycles generate a homology group of dimension two, i.e., when three minimal cycles are pairwise independent in homology, but threewise dependent. 
In this case, two generators are sufficient to span $H_1$ and, if their lengths are arranged pathologically, there is no principled way to choose two out of the three. \\
Suppose for example that three cycles $b_1,b_2$ and $b_3$ are such that
\[ \mu(b_1) < \mu(b_2) = \mu(b_3) \ \text{ and } \ [b_1] = [b_2] + [b_3] \]
In this case, both bases $\{ b_1, b_2 \}$ and $\{ b_1, b_3 \}$ span the same homology space, and are of equal minimal length. The minimality criterion fails in this case. \\
One could believe that such a configuration can only happen in the most general spaces, and that by imposing some mild hypotheses on the input data one could rule the pathology out. In fact the opposite is true, this degeneracy being possible even after enforcing very strong conditions on the data. \\
\noindent
\textbf{Counterexample} Even if $W$ is planar and an isometric embedding $W \hookrightarrow \mathbb{R}^2$ exists (i.e., the input planar weighted graph can be accurately drawn onto the plane), the minimal scaffold $\tilde{H}_{min}(W)$ needs not be unique. \\
In fact, consider complex $K$ arising from the geometric, planar graph in Fig. \ref{fig:PanelDraws}(d). Its homology $H_1(K)$ is generated by two cycles; possible generators are depicted in Fig. \ref{fig:PanelDraws}(e). Since the outer cycle $b_1$ is the shortest, and the two inner ones $b_2$ and $b_3$ are of equal length, the minimality criterion can not solve between $\{ b_1, b_2 \}$ and $\{b_1 , b_3 \}$, as both are acceptable minimal bases. 
The minimal scaffold (with or without draws) is not unique in this case. \\

\noindent
Clearly, the same could happen with more than three cycles, with a larger number of possibly ambiguous configuration. Therefore, if we allow for a high degree of symmetry in the input, this pathology could arise even in the rather tame context of planar graphs on $\mathbb{R}^2$. 
This issue is rather delicate, in the sense that not only the algorithm is unable to make a principled choice; it is not even capable of detecting when such a configuration takes place. 
In fact, this is more of a feature of homology than a flaw in the skeletonization framework: what our eyes see as different cycles are in fact homologically equivalent, and it is impossible to use homology to tell them apart. \\

\noindent
We however remark that, for complexes arising from real-world data, this type of configuration is actually pathological. 
Indeed, the following generality result holds \\

\noindent
\textbf{Proposition} Assume a point cloud $C = \{ X_i \}$ such that $X_i \sim U([0,1]^d)$ independently. Then, almost surely, the minimal scaffold $\mathcal{H}_{min}(W)$ (with or without draws) is unique. \\

\noindent
If the input point cloud is sampled uniformly at random in some $\mathbb{R}^ d$, then edge lengths are distributed according to an absolutely continuous probability law. Therefore, given two edges $e_1$ and $e_2$, $\mathbb{P}[ \mu(e_1) = \mu(e_2) ] = 0$. The same holds for any two non-identical cycles, and any two homology bases (being but finite sets of edges): the probability of them sharing the exact same length is zero. By finiteness of the input, at least one minimal homology basis exists and, by the above reasoning, almost surely this basis is unique for each filtration step. Then, with probability 1 the minimal scaffold is unique.  \\

\noindent
This result is actually quite general: whenever we can assume our input data to be subject to noise, then we are in principle allowed to rule out pathological same-length cycles. In these cases, the minimal scaffold is unique. \\

\noindent
We remark that this uniqueness result is compatible with the phenomenon of the concentration of measure: while for a very high-dimensional space or a very large number of points we know from theory that the distribution of length of edges concentrates towards its mean value, the probability of two edges (and hence two cycles) having the same length is still zero. One needs to be careful, however, that the probability of two cycles differing in length by less than some $\epsilon > 0$ could grow very rapidly with $\epsilon$.

\noindent
In summary, the minimal scaffold with draws $\tilde{\mathcal{H}}_{min}$ is well-defined up to some pathological circumstances, where it may depend on the ordering of the input. \\

\FloatBarrier

\section{Applications}
\label{sec:applications}

\noindent
As illustrative examples, we show here a few applications of the minimal scaffold. Through it, we obtain meaningful subsets of known networks in neuroscience, and rank their constituents by their ``topological importance". \\
    
\begin{figure}
        \centering
        \includegraphics[scale = 0.45]{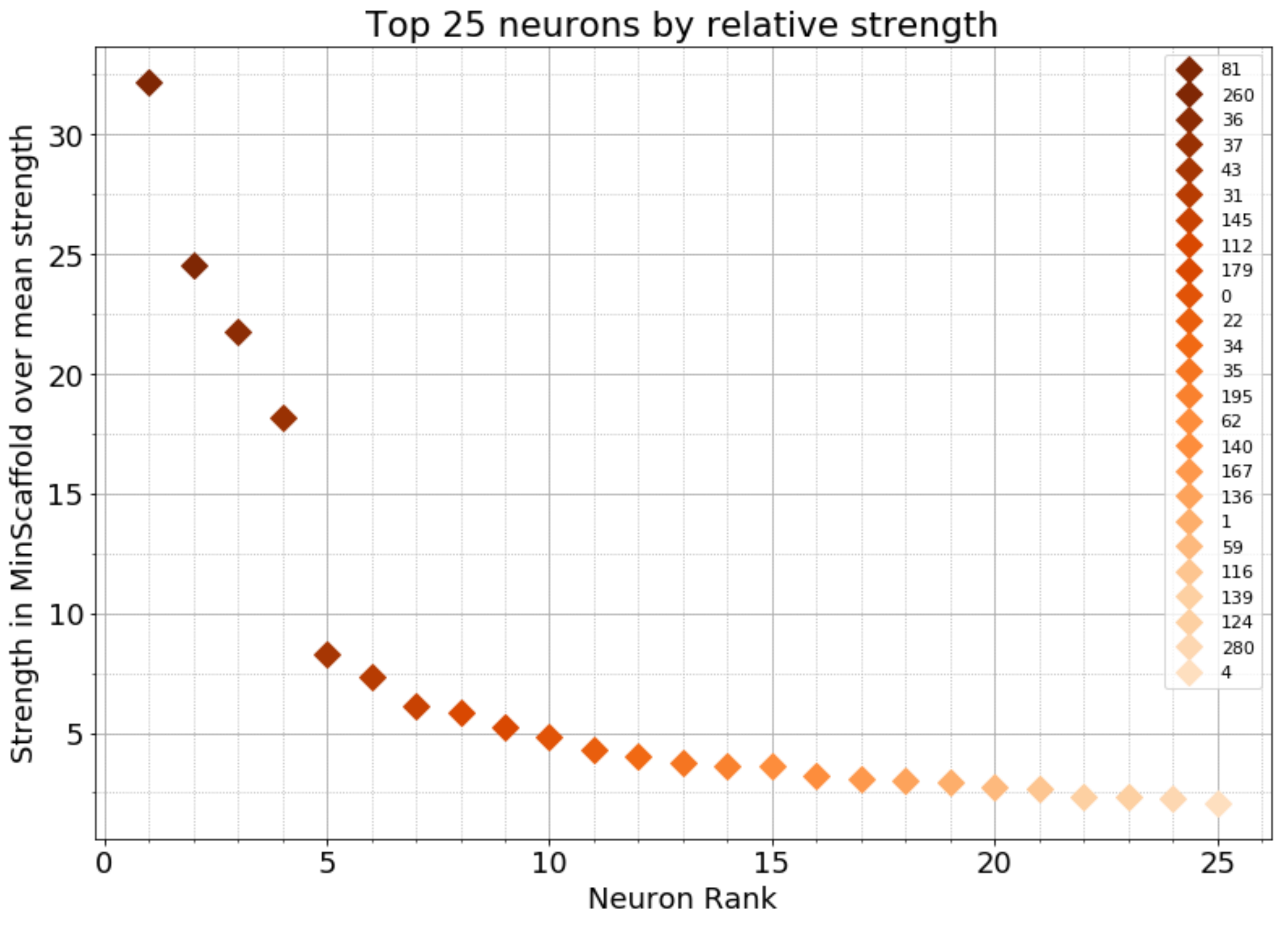}
        \caption{The top 25 neurons by relative node strength in the minimal scaffold over average strength in C. Elegans (mean $36.41$). Four neurons show a significantly higher relative strength than the others.}
        \label{fig:CETop}
    \end{figure}

\FloatBarrier
\noindent
The C. Elegans dataset is a correlation network of neural activations of the nematode worm Caenorhabditis Elegans. C. Elegans has become a model organism due to the unique characteristic of each individual sharing the exact same nervous system structure. \\
The input consists of a symmetric weighted adjacency matrix over 297 nodes, each representing a neuron. Edge weights represent (quantized) time correlations between the firing of neurons, ranging from 1 to 70. \\
The minimal homological scaffold of its brain map highlights the \emph{geometry} of the obstruction patterns, i.e., the precise areas where nervous stimuli are less likely to flow. We stress the improvement obtained by the minimal scaffold over the loose one, in that it is not only able to identify the \emph{presence} of a ``grey area" in the network, but it can as well provide a reliable boundary for it, and identify which neurons and inter-neuron links are responsible for information flowing around the obstruction. \\
As an interesting example, we see in Fig. \ref{fig:CETop} the top 25 neurons ranked in descending order of relative node strength (sum of weights of incident edges) with respect to the average node strength. We can identify four nodes, labeled 81, 260, 36, and 37, which hold a significantly higher relative strength than the rest. This implies their presence in many minimal cycles across several scales, hence suggesting that they play a crucial role in the fabric of information flow within the nematode's brain. \\

\begin{figure}
    \centering
    \includegraphics[scale=0.45]{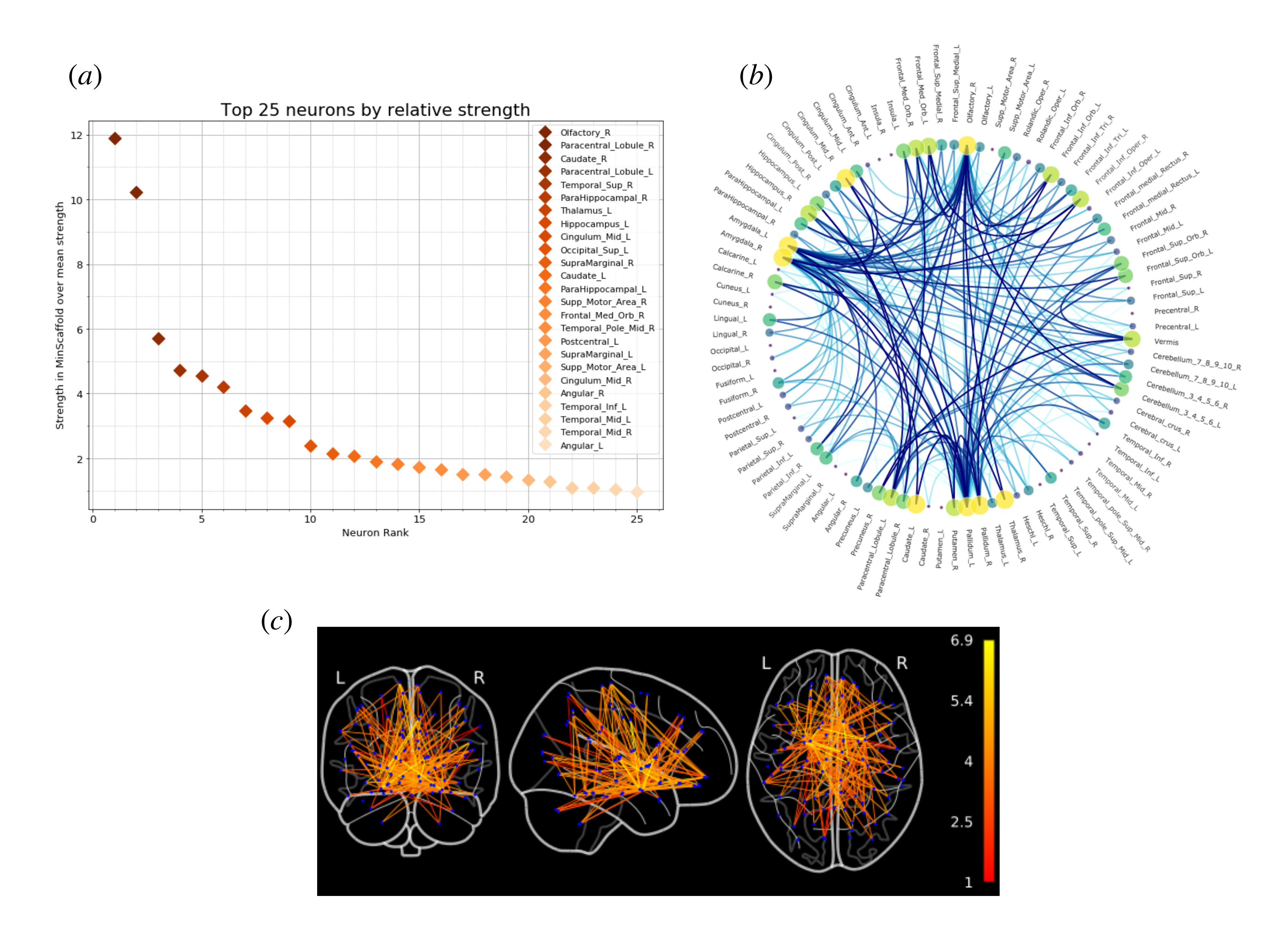}
    \caption{(a) The top 25 brain regions in the human brain by relative node strength in the minimal scaffold over average strength (mean $546.7$). Two neurons show significantly higher importance. (b) The chord diagram of the minimal scaffold. Node size represents node strength, edge color intensity represents weight in the scaffold. (c) The minimal scaffold embedded in the human brain, with regions accurately located, projected on the three coordinated planes. Edge color represents log-weight in the minimal scaffold (Log-scale for visualization purposes). }
    \label{fig:PanelBrain}
\end{figure}

\FloatBarrier
\noindent
The same type of analysis was repeated on the correlation network of brain activities in an 88-parcel atlas of the human brain, obtained through fMRI imaging at resting state. The data is courtesy of the Human Connectome Project (\cite{HumanConnectome}).  \\
Again, the minimal scaffold identifies which regions and links in the human brain are key bridges for the flow of information. Two parcels stand out (Fig. \ref{fig:PanelBrain}(a)) as particularly relevant for network topology. \\
For a relatively small network such as this, we can visualize the scaffold as a proper subnetwork by a chord diagram (Fig. \ref{fig:PanelBrain}(b)), with edge weight represented by color intensity and node strength by the size and color of the vertex. We stress that, starting from a virtually complete graph over 88 nodes, we reduce the size from 3828 edges to just 191, while preserving the topological structure.  \\ 
We can, as well, leverage libraries in computational neuroscience (\cite{NilearnPaper}) to embed the scaffold in the actual human brain, with regions correctly located, projected on the three coordinated planes. In Fig. \ref{fig:PanelBrain}(c), for visualization purposes color intensities represent log-weight in the scaffold. \\
To better highlight the value of the scaffold in signalling brain network function, we constructed a suitable null model of the functional network, as was done in \cite{Mastrandrea2017}. The technique consists in reshuffling the correlation matrix subject to the constraint of keeping a fixed spectrum, i.e. applying a random rotation, which guarantees the matrix remains positive semidefinite and hence a proper correlation matrix. An implementation of such a procedure can be found in \cite{Davies2000}. \\
The resulting randomized adjacency matrix is characterized by a vastly larger number of homological cycles than the original; so much so in fact that the computation of its minimal scaffold becomes cumbersome. However, even without computing them explicitly, we know for sure that the scaffolds of the original and randomized networks are totally different, specifically because they are built by aggregating two completely different persistence structures, i.e. the minimal scaffold does indeed highlight the functional information in the original dataset. \\

\noindent
The possible applications in which the minimal scaffold could provide novel insight into the structure of brain data are many: any relatively small correlation matrix could be either compressed or its patterns analyzed, as is often the case in EEG \cite{kuhlen2012EEG, ibanez2019topology,Tadic18Hyperbolicity, ibanez2019spectral} or neuronal \cite{GiustiPastalkovaCurto2015} studies, and in fMRI ones when using rather coarse atlases (e.g. \cite{termenon2016reliability,Tadic20Hubs}).
\FloatBarrier

\section{Comparison of Scaffolds}
\label{sec:comparison}

As the last contribution for this work, we consider a comparison between the minimal and loose scaffolds. \\
We have already pointed out that the minimal scaffold in general offers superior guarantees as a tool, both for network analysis and network skeletonization. On the other hand, the loose scaffold clearly has an advantage in terms of computational complexity: while it is in principle viable for most of the applications where persistent homology has been employed, the minimal scaffold, even adopting filtration-wise parallelization, requires a vastly larger amount of computational power, which effectively limits its range of application, unless run on dedicated, high-performance infrastructures. \\
A reasonable question to ask is the following. If one is interested not in the exact structure of the scaffold, but only in its statistical behaviour, could the loose scaffold provide a sufficient approximation of the minimal one? In a more concrete example, if instead of wondering exactly which nodes in a network are the most topologically important one is interested in the distribution of the degree sequence of the minimal scaffold, could the loose one come to one's help? \\
To answer this question, we have performed comparisons of several graph metrics in the two scaffolds of C. Elegans. Further, to gain insight into the general case, we have sampled two families of random graphs at different parameter values, one for geometric graphs (Random Geometric Graph), and one for non-geometric graphs (Weigthed Watts-Strogatz).

\subsection*{C. Elegans}
\FloatBarrier
For the C. Elegans dataset, we have compared the following graph metrics of the minimal and loose scaffolds:
\begin{enumerate}
    \item Degree Sequence
    \item Node Strength
    \item Betweeness Centrality
    \item Closeness Centrality
    \item Eigenvector Centrality
    \item Clustering Coefficients
    \item Edge weights
\end{enumerate}
    
\begin{figure}
    \centering
    \includegraphics[scale=0.43]{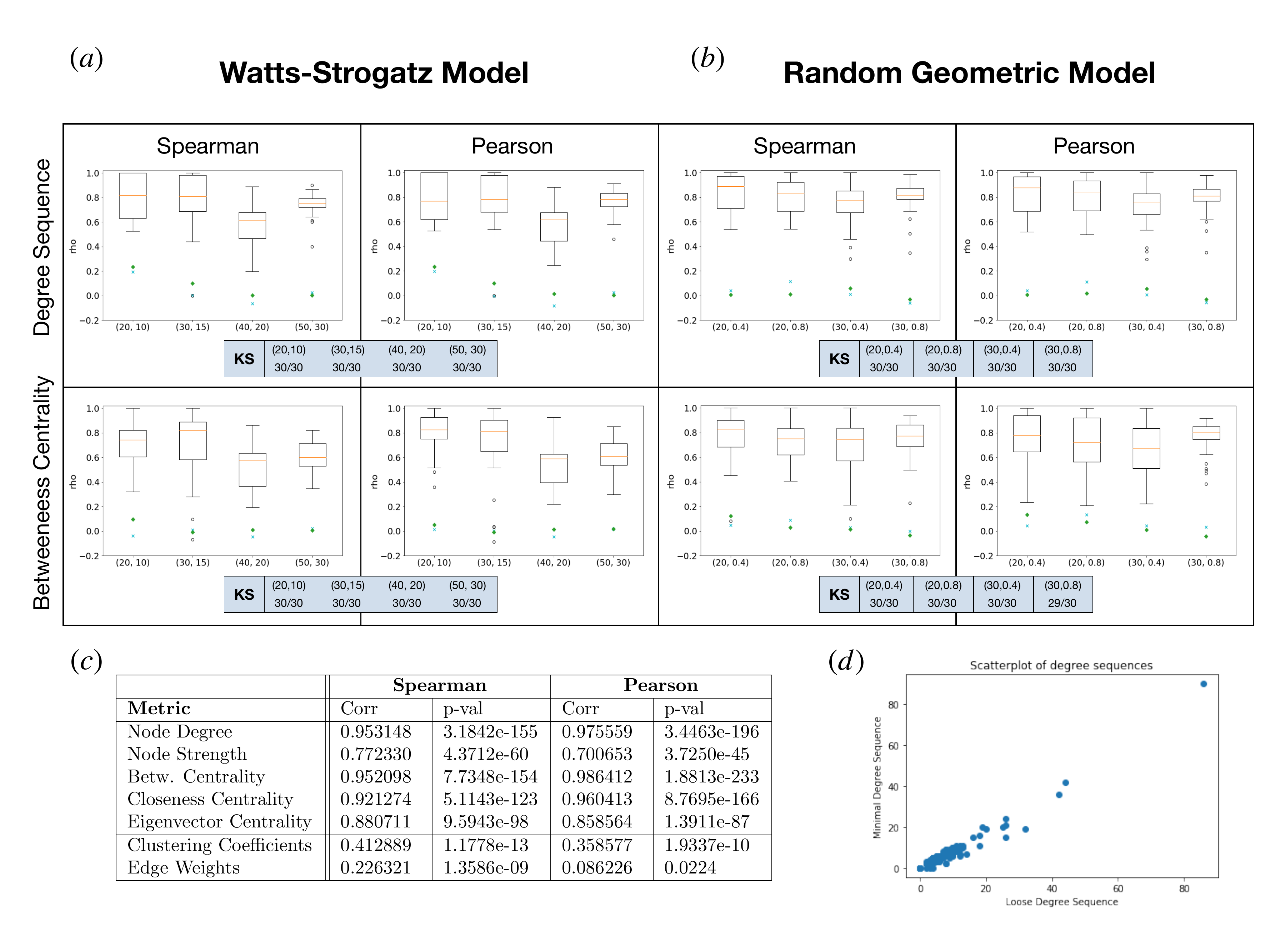}
    \caption{Correlations between the minimal and loose scaffold. (a) Comparison in the weighted Watts-Strogatz model. Degree sequence and betweenness centrality in the two scaffolds are compared, using Pearson and Spearman correlation coefficients. Each box is computed over a sample of 30 weighted Watts-Strogatz random graphs, with parameters as reported on the x-axis: the pair $(N,k)$ indicates a WS model on $N$ nodes, with $k$ stubs to rewire. The rewiring probability is $0.025$. The cyan crosses and the green diamonds represent the average correlation value against the loose and minimal null models, respectively. (b) Comparison in the random geometric model. Again, Pearson and Spearman correlation coefficients of the degree sequence and betweenness centrality in the two scaffolds are compared. Each box is computed over a sample of 30 random geometric graphs, with parameters as reported on the x-axis: the pair $(N,t)$ indicates a graph on $N$ nodes sampled uniformly at random in the $[0,1]^2$ square. $t$ is the connectivity distance threshold. The cyan x's and the green diamonds represent the average correlation value against the loose and minimal null models, respectively. The darker boxes in panels (a) and (b) report, for their respective model and for each metric and parameter values, the fraction of the sampled instances for which the Kolmogorov-Smirnov test was inconclusive ($p$ value > $0.05$). (c) Correlation tests for several network metrics show significant capabilities of the standard scaffold to reproduce certain statistical properties of the minimal one in C. Elegans. At the same time, due to different construction mechanisms, others are unreliable. (d) Scatterplot of the degree sequence of neurons of C. Elegans in the minimal scaffold versus in the loose one.}
    \label{fig:PanelBoxplots}
\end{figure}

\noindent
Results (reported in the Table of Fig. \ref{fig:PanelBoxplots}(c)) indicate that, for metrics 1 to 5, the two scaffolds are very well correlated. So for example the cheap, loose scaffold is a reliable proxy of the distribution of the ``true" degree sequence (scatterplot in Fig. \ref{fig:PanelBoxplots}(d)). \\ 
We instead observe poor correlation of edge weights and clustering coefficients. The first one is not unexpected, since the edge weighting procedure is conceptually different in the two scaffolds: while in the minimal one we consider a different basis for each filtration step, the loose scaffold considers bases of the persistent homology space, drastically reducing the number of cycles considered. To make it clearer, in general set $B^*$ has cardinality much larger than the dimension of $PH_1$. It is therefore explicable that the distributions of edge weights do not generally agree. \\
Clustering coefficients, on the other hand, are a measure of how ``triangular" a graph is around a given node. As remarked in Section \ref{sec:minimal-scaffold}, another consequence of assembling the scaffold from the minimal bases of the $H_1$'s is that a large number of artificial triangles appear around cycles. In this case too, therefore, the poor correlation is easily explained.

\subsection*{Random Graphs}
Drawing inspiration from \cite{SizBassClass17}, we repeat the analysis on random graph samples. \cite{SizBassClass17} divides random networks into two categories: those created from edge weighting schemes and those created from points in the Euclidean space. We have chosen to analyze the weighted Watts-Strogatz (WS) model as representative of the first class, and the geometric random model as representative of the second. We remark that weighting needs to be introduced in order to compute persistence; while for geometric graphs this simply requires computing the Euclidean distance, for the Watts-Strogatz model it requires an ad-hoc procedure that is described in detail in the supplemental material of \cite{SizBassClass17}. \\
We briefly recall that a WS graph is parametrized by the number of nodes, by the number of stubs to rewire, and by the rewiring probability. A random geometric graph is instead parametrized by the number of points to sample (uniformly) in $[0,1]^d$, and by a cutoff value that acts as distance threshold, beyond which no edge is introduced. \\
In both cases, we observe good agreement on key statistics, as reported in Fig. \ref{fig:PanelBoxplots}(a) and (b). Each bar is obtained by computing the correlation of the reported statistic on a sample of 30 random graphs of the reported model, with parameters as indicated on the x-axis.\\
For comparison, two null models are built for each instance of the minimal and loose scaffolds in the sample, by constructing an Erd\H{o}s-R\'{e}nyi random graph on the same vertex set, one with the same number of edges as the minimal scaffold, and one with the same number as the loose one. The correlation is computed of each statistic between the minimal scaffold and the loose null model and between the loose scaffold and the minimal null model. The average of these correlations is reported on the boxplots to act as a baseline value, highlighting that the two scaffolding procedures agree with each other by more than just statistical noise. \\
For a finer analysis, we have performed a two-sample Kolmogorov-Smirnov test comparing the distribution of the given metrics in the minimal and loose scaffolds, for all parameter values of the two random models. We consider the Kolmogorov-Smirnov test to be inconclusive if its $p$ value exceeds a threshold of $0.05$, in which case one cannot confidently reject the null hypothesis that the samples are drawn from the same distribution. In Fig. \ref{fig:PanelBoxplots} panels (a) and (b), the darker boxes report for each parameter choice and metric the fraction of samples for which the test was inconclusive: in all cases except one, the KS test could not distinguish between the distribution of the graph statistic between the minimal and loose scaffolds, strengthening the indication of a good agreement between the two.

\FloatBarrier

\subsection*{nPSO Random Graph Model}
\begin{figure}
    \centering
    \includegraphics[scale=0.43]{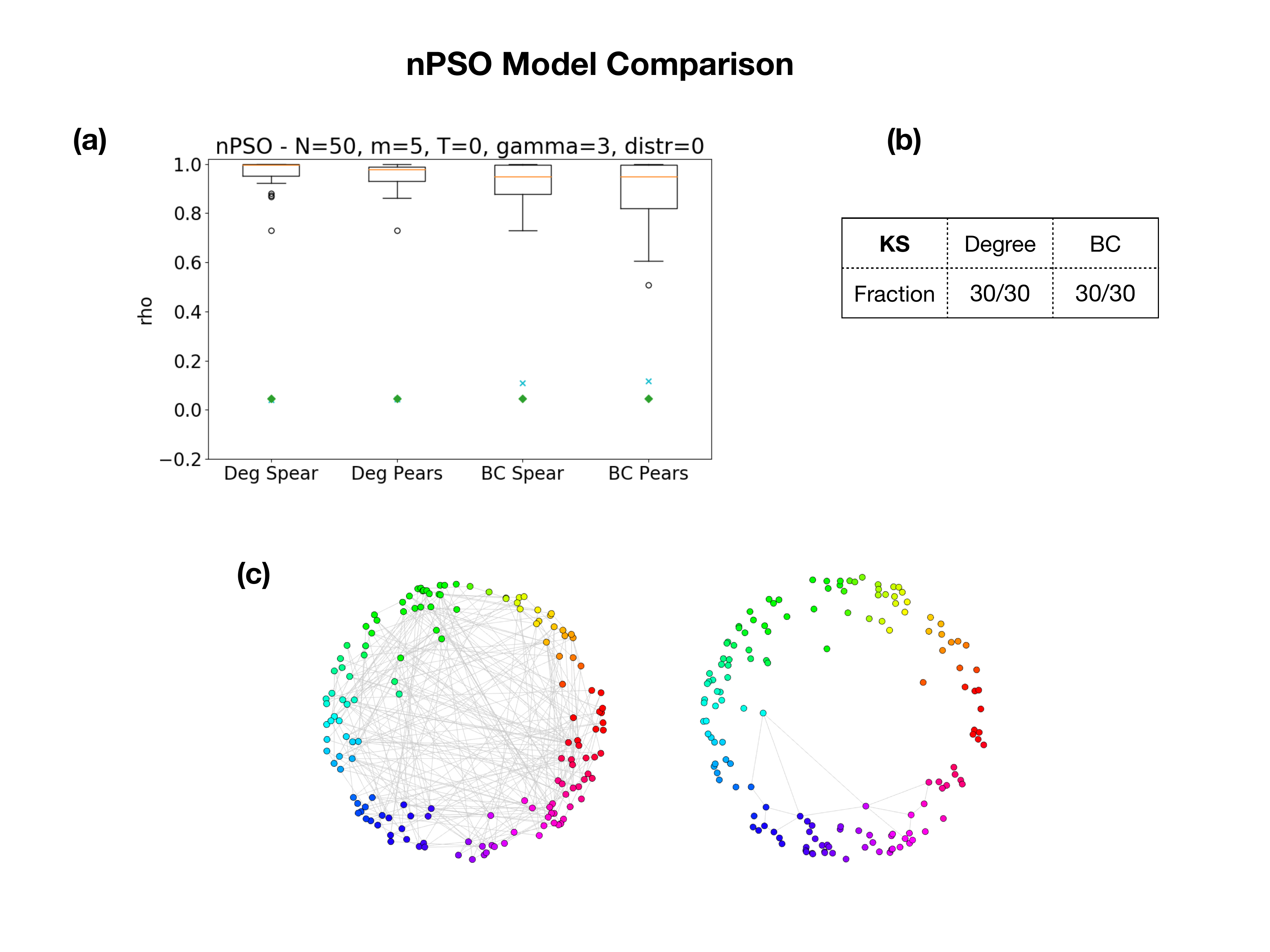}
    \caption{Comparison of the minimal and loose scaffold for nPSO random model. (a) Degree sequence and betweenness centrality in the two scaffolds are compared using Pearson and Spearman correlation coefficients. Each box is computed over a sample of 30 nPSO instances, with the following parameters: 50 nodes, average degree 10 ($m=5$), 0 temperature, power-law exponent $\gamma=3$, and uniform distribution of angular coordinates. The cyan crosses and green diamonds represent the average correlations against the loose and minimal null models respectively, as in Fig. \ref{fig:PanelBoxplots}. In panel (b), the table reports, for the degree and betweenness centrality distributions, the fraction of Kolmogorov-Smirnov test that could not reject the hypothesis of the two samples coming from the same distribution. This has always been the case for each sampled instance and both metrics. (c) A graphical depiction of an instance of the nPSO model with parameters $N=150, m=2,T=5,\gamma=3$ and uniform distribution on the left. On the right, the corresponding minimal scaffold.}
    \label{fig:PanelNPSO}
\end{figure}

A modern random graph model, which has recently gained traction in network science for its ability to concurrently tune several parameters of interest in modeling real networks, is the Nonuniform Popularity-Similarity model. Introduced in \cite{Muscoloni2018Nonuniform}, it builds upon a sequence of increasingly refined generative models to provide all the key structural properties of real-world graphs, such as scale-freeness,  small-worldness and community structure. We therefore set out to employ it as benchmark in our comparison of the minimal and loose scaffolds. \\
In general, networks which display hyperbolic geometries tend to have a rather tree-like structure, with a certain scarcity of cycles. It is straightforward that, in the absence of a significant structure of persistent homology, the loose and minimal scaffolds will agree to high degree for at least two reasons: the low number of cycles forces the loose scaffold to localize onto the few available holes, hence resembling the minimal, and secondly the scarcity of homology makes for a comparison between two mostly empty sets. \\
Following the lead of \cite{Muscoloni2018Leveraging}, we tuned the nPSO model parameters in order to empirically maximize the persistent homology structure, so as to make the comparison the most significant possible. As reported in Fig. \ref{fig:PanelNPSO}, we observe again good ability of the scaffolds to proxy each other across the metrics analyzed, significantly higher than with respect to a null model, for a sample with parameters $N=50, m=2,T=5,\gamma=3$ and uniform distribution. A Kolmogorov-Smirnov test was also performed, as in the previous section, where a $p$-value higher that $0.05$ indicates that the distribution of degrees and betweeness centralities in the minimal and loose scaffold cannot be confidently distinguished. This was the case for all the samples we tested.

\FloatBarrier

\section{Conclusions}
\label{sec:conclusions}
We provided a new method of network analysis and skeletonization, based on the computation of minimal homology bases. This new construction fills a significant gap in previous literature, in that it yields, in all but some pathological cases, a well-defined and unique subgraph, acting as a reasonable ground truth for comparison with the  previous construction. It can be employed in a range of applications, both to identify crucial and weak links in a network, and to obtain compressed and topologically sound representations of the input. 
It also allows to evaluate the reliability of other scaffolding procedures with respect to said ground truth: we have observed that, for some applications, the loose scaffold can be deemed a sufficiently accurate tool, while not incurring in as cumbersome a computational load. \\
We foresee that the subject of homological skeletonization is not yet concluded. Other approaches to finding canonical generators of homology are possible (for example in \cite{Kurlin1D2015} and \cite{boissonnat2018}), and we plan to investigate them further in subsequent works. \\
A question which remains open and could be worthy of further work is the following: could one construct a sensible "entropy" functional on the space of cycles, so as to obtain a strictly unique, minimally-represented basis that is in the most likely?

\bibliography{bibliography.bib}

\section*{Availability of data and material}
The C. Elegans dataset analysed during the current study is available and included in the GitHub repository MinScaffold, \href{https://github.com/marcoguerra192/MinScaffold}{https://github.com/marcoguerra192/MinScaffold}. \\
The Human Connectome Project dataset is available from the page \\ \href{http://www.gipsa-lab.grenoble-inp.fr/~sophie.achard/Brain_connectivity_network}{http://www.gipsa-lab.grenoble-inp.fr/~sophie.achard/Brain\_connectivity\_network}

\section*{Acknowledgements}

MG, ADG, UF, and FV acknowledge the support from the Italian MIUR Award ``Dipartimento di Eccellenza 2018-2022" - CUP: E11G18000350001 and the SmartData@PoliTO center for Big Data and Machine Learning.
GP acknowledges partial support from Intesa Sanpaolo Innovation Center. The funder had no role in study design, data collection, and analysis, decision to publish, or preparation of the manuscript. \\

\noindent
The authors acknowledge Iacopo Iacopini for kindly sharing a Python library for plotting simplicial complexes, available on GitHub ( github.com/iaciac/py-draw-simplicial-complex). We further acknowledge the python library Nilearn (\cite{NilearnPaper}) for the brain image visualization code. 
We would also like to thank Paola Siri for useful discussions.

\section*{Author contributions statement}

MG, ADG, UF, GP, and FV conceived and designed the study, performed the analysis and wrote the manuscript. All authors read and approved the final manuscript.

\section*{Additional information}

\textbf{Competing interests} The authors declare that they have no competing interests.

\end{document}